\documentclass[english]{article}
\usepackage[T1]{fontenc}
\usepackage[latin9]{inputenc}
\usepackage{ifsym}
\usepackage{amsmath}
\usepackage{amssymb}
\usepackage{esint}

\usepackage{babel}

\begin{document}

\title{{\normalsize Semigroups of finite-dimensional random projections}}

\author{{\normalsize Andrey A. Dorogovtsev}}

\maketitle
Institute of Mathematics Ukrainian Academy of Sciences

adoro@imath.kiev.ua
\begin{abstract}
{\normalsize In this paper we present a complete description of a
stochastic semigroup of finite-dimensional projections in Hilbert
space. The geometry of such semigroups is characterized by the asymptotic
behavior of the widths}\emph{\normalsize{} }{\normalsize of compact
subsets with respect to the subspaces generated by the semigroup operators. }{\normalsize \par}
\end{abstract}

\section{{\normalsize Introduction.}}

This paper is devoted to a very special type of stochastic operator-valued
semigroups, namely to semigroups of finite-dimensional projections.
The study of stochastic operator-valued semigroups was originated
from the works of A.V. Skorokhod {[}1,2], where he gave a representation
of such a semigroup as a solution to a stochastic differential equation
with an operator-valued martingale or with a process with independent
increments. Skorokhod treated stochastic operator-valued semigroup
as the Dolean exponent for some operator-valued martingale. This approach
requires that the mean operators have continuous inverse. Semigroups
of finite-dimensional projections, which we study in this paper, do
not satisfy this condition. Moreover, the nature of semigroups of
finite-dimensional projections is different from those studied in
{[}1,2]. They arose as operators describing shifts of functions or
measures along a stochastic flow. In {[}2] A.V.Skorokhod briefly noted
that in the case when a stochastic flow is generated by a stochastic
differential equation, such semigroups must satisfy a stochastic differential
equation with unbounded operators in coefficients, but detail investigation
was not provided.

Stochastic flows with coalescence in general can not be represented
as a solution to a stochastic differential equation. An example of
such flows is the Arratia flow of Brownian particles {[}4]. It is
a partial case of so-called Harris flows {[}5], which consist of Brownian
particles with spatial correlation depending on the difference between
the positions of the particles. In contrast to the flows generated
by stochastic differential equations such flows can lose the gomeomorphic
property {[}5]. Moreover, the Arratia flow maps every bounded interval
into a finite number of points. Consequently, for the investigation
of coaleshing stochastic flows we can not apply the tools from differential
geometry used for smooth stochastic flows (about smooth stochastic
flows see {[}6]). 

One of the possible approaches to the understanding of the geometry
of stochastic flows with coalescence is the investigation of random
operators describing shifts of functions or measures along a flow.
It will be shown in the next section that such random operators can
have a finite-dimensional image almost surely and be unbounded. The
corresponding semigroup can have a complicated structure. So in this
article we propose to consider a stochastic semigroup consisting from
finite-dimensional projections. It occurs that such a semigroup in
Hilbert space can be described in terms of a certain random Poisson
measure on the product of the positive half-line and some orthonormal
basis (Theorem 3.1). The geometry of such semigroup can be characterized
in terms of widths of compact sets with respect to it. The asymptotic
behavior of such widths for two types of compacts we establish in
Section 4. The study of the semigroup related to a stochastic flow
with coalescence in the general case is the subject to the ongoing
work.

\section{{\normalsize Random projections in Hilbert space.}}

Let us recall the definition of the random operator. Let $H$ be a
real separable Hilbert space with an inner product $(.,.)$ and $L_{2}(\Omega,P,H)$
be the space of square-integrable random elements in $H$.

\textbf{\large Definition 2.1.} Strong random operator in $H$ is
a continuous linear map from $H$ to $L_{2}(\Omega,P,H)$ .

In this section we will consider a special case of the strong random
operators. Namely we will consider finite-dimensional random projections.
They are bounded random operators. 

\textbf{\large Definition 2.2.} A strong random operator $A$ is a
bounded random operator if there exists a family $\{A_{\omega},\omega\in\Omega\}$
of deterministic bounded linear operators in $H$ such, that
\begin{description}
\item [{$\forall x\in H:\:\:(Ax)_{\omega}=A_{\omega}x$}] $a.s.$
\end{description}
Note that boundedness of a strong random operator and finiteness of
its image dimension are not connected, as it is demonstrated in the
next example. 

\textbf{\large Example 2.1.} Suppose that $H=L_{2}([0;1]),\:\theta$
is a random variable uniformly distributed on $[0;1].$ Define a random
linear functional on $H$ as follows \[
H\ni f\longmapsto\Phi(f)=f(\theta).\]
Then \[
E\Phi(f)^{2}=\int_{0}^{1}f(s)^{2}ds.\]
Consequently, $\Phi$ is continuous in the square-mean. Let us check
that $\Phi$ is not a random bounded functional. Suppose, that it
is and denote by $\Phi_{\omega},\:\omega\in\Omega$ the corresponding
family of continuous linear functionals on $H$. Let us identify the
functional $\Phi_{\omega}$ with the function from $H$ for every
$\omega\in\Omega$ . Let $\{f_{n}\in H$, $n\geq1\}$ be a sequence
of i.i.d. random variables with a finite second moment when considered
on the standard probability space $[0;1]$ with the Lebesgue measure.
Suppose also that \[
\int_{0}^{1}f_{1}(s)^{2}ds=1,\:\int_{0}^{1}f_{1}(s)ds=0.\]
Then the sequence $\{\Phi(f_{n}),\: n\geq1\}$ is the sequence of
i.i.d. random variables. From the other side \[
\Phi(f_{n})_{\omega}=(f_{n},\Phi_{\omega})\rightarrow0,\: n\rightarrow\infty.\]
This contradiction proves our statement.

In what follows we say that $A$ is a random Hilbert-Schmidt or nuclear
or finite-dimensional operator if the corresponding family $\{A_{\omega},\omega\in\Omega\}$
consists of Hilbert-Schmidt, nuclear or finite-dimensional operators.
Let us start with a characterization of random Hilbert-Schmidt operators. 

\textbf{\large Theorem 2.1.} A strong random operator $A$ is a random
Hilbert-Schmidt operator if and only if for some orthonormal basis
$\{e_{n};n\geq1\}$ in $H$ 

\begin{equation}
\sum_{n=1}^{\infty}||Ae_{n}||^{2}<\infty.\end{equation}

\textbf{\large Proof. }If $A$ is a random Hilbert-Schmidt operator,
then the condition of the theorem holds obviously. Suppose that strong
random operator $A$ satisfies (1). Define a subset $\Omega_{0}$
of the probability space as follows
\begin{description}
\item [{$\forall\omega\in\Omega_{0},\alpha_{1},...,\alpha_{n}\in R:$}] \[
\sum_{k=1}^{n}\alpha_{k}Ae_{k}(\omega)=A(\sum_{k=1}^{n}\alpha_{k}e_{k})(\omega),\:\sum_{n=1}^{\infty}||Ae_{n}||^{2}(\omega)<\infty.\]

\end{description}
Now for every $\omega\in\Omega_{0}$ define a Hilbert-Schmidt operator
$A_{\omega}$ as follows\[
{\displaystyle A_{\omega}x=\sum_{n=1}^{\infty}(x,e_{n})}(Ae_{n})(\omega).\]
For $\omega\notin\Omega_{0}$ put $A_{\omega}=0$. Then $A$ satisfies
Definition 2.2 with the family $\{A_{\omega},\omega\in\Omega\}.$
\textifsymbol[ifgeo]{48}

The previous statement allows us to characterize random finite-dimensional
projections. 

\textbf{\large Theorem 2.2.} Suppose that a strong random operator
$A$ satisfies the condition of Theorem 2.1 and the following conditions
hold:
\begin{itemize}
\item for all $x,y\in H\,\quad(Ax,y)=(Ay,x),\,(Ax,x)\geq0,$
\item $A=A^{2}.$
\end{itemize}
Then $A$ is a random finite-dimensional projection.

\textbf{\large Proof.} The proof follows immediately from the well-known
characterization of projections in a Hilbert space and the fact that
any Hilbert-Schmidt projection is finite-dimensional. \textifsymbol[ifgeo]{48}

Let us consider some examples of random finite-dimensional projections.

\textbf{\large Example 2.2.} Let $n(t),\: t\in[0;1]$ be a Poisson
process. Denote by $\tau_{1},...,\tau_{\nu}$ the subsequent jumps
of $n.$ The intervals $[0;\tau_{1}),\ldots,[\tau_{\nu};1]$ generate
a finite $\sigma-$field $\mathit{\mathcal{A}}$. Define a random
operator $A$ in $L_{2}([0;1])$ as a conditional expectation with
respect to the $\sigma-$field $\mathit{\mathcal{A}}$. Then $A$
is a random finite-dimensional projection.

\textbf{\large Example 2.3.} Let $(\mathit{X},\rho)$ be a Polish
space and $\mu$ be a probability measure on the Borel $\sigma-$field
in $\mathit{X}$. Consider a measurable map $\phi:\mathit{X}\times\Omega\mapsto\mathit{X}$
such, that for every $\omega\in\Omega$ the image $\phi(\mathit{X},\omega)$
contains a finite number of elements. Define a random operator $A$
in $L_{2}(\mathit{X},\mu)$ as a conditional expectation with respect
to the $\sigma-$field generated by $\phi$. Then $A$ is a finite-dimensional
random projection related to $\phi.$ 

\textbf{\large Remark.} If the map $\phi$ has not the finite image
the operator $A$ still will be well-defined random projection but
not finite-dimensional.

\section{{\normalsize Semigroups of projections.}}

In this section we introduce the notion of semigroups of random projections
which is the main object of investigation in the article.

\textbf{\large Definition 3.1.} A family of random bounded operators
$\{G_{s,t},\:0\leq s\leq t<\infty\}$ is referred to as a semigroup
if the following conditions hold:
\begin{enumerate}
\item For any $s,t,r\geq0:\;$ $G_{s,t}$ and $G_{s+r,t+r}$ are equidistributed.
\item For any $x\in H:\;$ ${\displaystyle E||G_{0,t}x-x||^{2}\mapsto0,\: n\mapsto\infty.}$
\item For any $0\leq s_{1}\leq\cdots\leq s_{n}<\infty:\; G_{s_{1},s_{2}}$
,$\ldots$ ,$G_{s_{n-1},s_{n}}$ are independent.
\item For any $r\leq s\leq t:\;$ $G_{r,s}G_{s,t}=G_{r,t}$, $G_{r,r}=I$,
where $I$ is an identity operator.
\end{enumerate}
To illustrate a connection of the stochastic semigroup with stochastic
flows consider the following example. Suppose, that $X$ is a Polish
space. Denote by $\mathbf{M}$ the space of all finite signed measures
on the Borel $\sigma-$field in $X$ equipped with the topology of
weak convergence. $\mathbf{M}$ is a linear topological space. 

\textbf{\large Definition 3.2.} A family of measurable mappings $\phi_{s,t}:X\times\Omega\mapsto X$
, $0\leq s\leq t<\infty$ is referred to as a random flow on $X$
if the following conditions hold:
\begin{itemize}
\item For any $0\leq s_{1}\leq s_{2}\leq\ldots s_{n}<\infty:\;\phi_{s_{1},s_{2}}$
,$\ldots$ ,$\phi_{s_{n-1},s_{n}}$ are independent.
\item For any $s,t,r\geq0:\;$$\phi_{s,t}$ and $\phi_{s+r,t+r}$ are equidistributed.
\item For any $r\leq s\leq t$ and $u\in X:\;$ $\phi_{r,s}\phi_{s,t}(u)=\phi_{r,t}(u)$,
$\phi_{r,r}$ is an identity map.
\item For any $u\in X:\;$$\phi_{0,t}(u)\mapsto u$ in probability when
$t\mapsto0.$
\end{itemize}
\textbf{\large Example 3.1.} Suppose that the random flow $\{\phi_{s,t},\:0\leq s\leq t<\infty\}$
on $X$ has the following additional property. For any $t$ the function
$\phi_{0,t}$ is continuous with probability one. Define the operator
$G_{s,t}$ in $\mathbf{M}$ by the formula\[
G_{s,t}(\mu)=\mu\phi_{s,t}^{-1}.\]
It can be easily checked that the family $\{G_{s,t},\:0\leq s\leq t<\infty\}$
satisfies the analog of Definition 3.1 for a linear topological space. 

In the next example we consider a semigroup of random finite-dimensional
projections in Hilbert space. 

\textbf{\large Example 3.2.} Let $H$ be Hilbert space with an orthonormal
basis $\{e_{k},\: k\geq1\}.$ Consider the sequence $\{n_{k},\: k\geq1\}$
of independent Poisson processes with intensities $\{\lambda_{k},\: k\geq1\}.$
Suppose that\begin{equation}
\forall\rho>0:\quad{\displaystyle \sum_{k=1}^{\infty}\exp(-\rho\lambda_{k})<+\infty}\end{equation}
 Define for every $k\geq1$ and $0\leq s\leq t$ \[
{\displaystyle \nu_{s,t}^{k}=\begin{cases}
0,\: n_{k}(t)-n_{k}(s)>0,\\
1,\: n_{k}(t)-n_{k}(s)=0.\end{cases}}\]
Finally define the projection $G_{s,t}$ as follows\[
{\displaystyle G_{s,t}(u)=\sum_{k=1}^{\infty}(u,e_{k})\nu_{s,t}^{k}e_{k}.}\]

Condition (2) implies that $G_{s,t}$ is a finite-dimensional projection
with probability one. The conditions of Definition 3.1 trivially hold. 

The next lemma shows that deterministic semigroup of finite-dimensional
projections does not exist.

\textbf{\large Lemma 3.1.} Suppose that $\{G_{t},\:0\leq t<\infty\}$
is a strongly continuous semigroup of bounded operators in separable
Banach space $\mathcal{B}$. Assume that $\dim G_{t}(\mathcal{B})<\infty$
for every $t>0$. Then $\dim\mathcal{B}<\infty.$ 

\textbf{\large Proof.} Define the function $\nu(t)=\dim G_{t}(\mathcal{B}),\: t>0.$
It is clear that this function is decreasing, takes integer values
and\[
{\displaystyle \lim_{t\rightarrow0}\nu(t)=+\infty.}\]
Let $t_{0}$ be a positive point of jump for the function $\nu$.
Then there exists a nonzero element $x\in G_{t_{0}}(\mathcal{B})$
such that $x\notin G_{t}(\mathcal{B})$ for all $t>t_{0}.$ Since
$\{G_{t},\:0\leq t<\infty\}$ is a semigroup, then for arbitrary $s>0$
$G_{s}(x)=0.$ This contradicts to the strong continuity of $G$.
\textifsymbol[ifgeo]{48}

Actually Example 3.2 describes the unique possibility of the construction
of the semigroup of random finite-dimensional projections in Hilbert
space. 

\textbf{\large Theorem 3.1.} Let $\{G_{s,t},\:0\leq s\leq t<\infty\}$
be a semigroup of random finite-dimensional projections in separable
Hilbert space H. Then there exists an orthonormal basis $\{e_{k},\: k\geq1\}$
in $H$ and Poisson processes $\{n_{k},\: k\geq1\}$ which have jointly
independent increments, such that \[
{\displaystyle {\displaystyle G_{s,t}(u)=\sum_{k=1}^{\infty}(u,e_{k})\nu_{s,t}^{k}e_{k},}}\]
where for every $k$ $\nu_{s,t}^{k}$ is built from $n_{k}$ exactly
in the same way as in Example 3.2. 

\textbf{\large Proof.} Consider two projections $R_{1},\: R_{2}$
in $H$ such that their product $Q=R_{1}R_{2}$ is a projection. Then
$R_{2}R_{1}=R_{1}R_{2}$. To check this relation introduce the notations
$R_{i}(H)=L_{i},\: i=1,2.$ Suppose, that $u\in H$ is such that $||u||=||$$R_{1}R_{2}u||$.
Then $||u||=||R_{2}u||,\:||u||=||R_{1}u||$ i.e. $u\in L_{1}\cap L_{2}$.
This means that $Q$ is a projection on $L_{1}\cap L_{2}$. Since
the subspace $L_{2}$ can be represented as $L_{2}=L_{2}^{'}\oplus L_{1}\cap L_{2}$,
then $R_{1}(L_{2}^{'})=\{0\}$. This imply that $R_{2}R_{1}=R_{1}R_{2}$.
Now consider the semigroup $\{G_{s,t},\:0\leq s\leq t<\infty\}$.
From the definition of the semigroup and the previous considerations
one can conclude that for every $s_{1}\leq s_{2}\leq s_{3}$ the operators
$G_{s_{1},s_{2}}$ and $G_{s_{2},s_{3}}$ commute with probability
one. Since the space $H$ is separable, then there exists a sequence
$\{\triangle_{kn};\: k,\: n\geq1\}$ of increasing partitions of the
probability space which generates $\sigma(G_{s_{1},s_{2}})$ . The
following integrals define bounded operators in $H$\[
Q_{s_{1},s_{2}}^{kn}u=\frac{1}{P(\triangle_{kn})}\int_{\triangle_{kn}}G_{s_{1},s_{2}}uP(d\omega),\: u\in H.\]
The norm of $Q_{s_{1},s_{2}}^{kn}$ is less or equal to one. Consider
the following random operators \[
G_{s_{1},s_{2}}^{n}=\sum_{k=1}^{\infty}Q_{s_{1},s_{2}}^{kn}1_{\triangle_{kn}},\: n\geq1.\]
 It can be easily checked that for every $u\in H$ with probability
one \[
G_{s_{1},s_{2}}^{n}u\rightarrow G_{s_{1},s_{2}}u,\: n\rightarrow\infty.\]
The random operators $G_{s_{1},s_{2}}^{n}$ take commuting values.
Namely, for the arbitrary $k,\: n,\: s_{1},\: s_{2}$ and $l,\: m,\: t_{1},\: t_{2}$
we have \[
Q_{s_{1},s_{2}}^{kn}Q_{t_{1},t_{2}}^{lm}=Q_{t_{1},t_{2}}^{lm}Q_{s_{1},s_{2}}^{kn}.\]
 To prove this relation note that the values of $Q_{s_{1},s_{2}}^{kn}$
depend only from the distribution of $G_{s_{1},s_{2}}$. Consequently,
it is enough to consider the case, when $0=s_{1}<s_{2}=t_{1}<t_{2}=s_{2}-s_{1}+t_{2}-t_{1}$.
Now the subsets which we use in the definition of $Q_{s_{1},s_{2}}^{kn}$
and $Q_{t_{1},t_{2}}^{lm}$ are independent. Hence \[
\int_{\triangle_{kn}}G_{s_{1},s_{2}}P(d\omega)\int_{\triangle_{lm}}G_{t_{1},t_{2}}P(d\omega)=\]
\[
=\int_{\Omega}G_{s_{1},s_{2}}G_{t_{1},t_{2}}1_{\triangle_{kn}}1_{\triangle_{lm}}P(d\omega)G_{s_{1},s_{2}}P(d\omega)=\]
 \[
=\int_{\Omega}G_{t_{1},t_{2}}G_{s_{1},s_{2}}1_{\triangle_{lm}}1_{\triangle_{kn}}P(d\omega)=\int_{\triangle_{lm}}G_{t_{1},t_{2}}P(d\omega)\int_{\triangle_{kn}}G_{s_{1},s_{2}}P(d\omega).\]
Also note that the operators $Q_{s_{1},s_{2}}^{kn}$ are self-adjoint
and non-negative. Consequently, one can build a countable family $\Gamma$
of commuting self-adjoint non-negative operators such that random
projections $G_{s_{1},s_{2}}$ can be approximated by random operators
with the values from $\Gamma$. Let us verify that we can choose the
family $\Gamma$ in such a way that it consists of the nuclear operators.
Since our projections are finite-dimensional, then for arbitrary $s_{1}\leq s_{2}$
\[
trG_{s_{1},s_{2}}<+\infty.\]
Then truncating the sets $\triangle_{kn}$ to their intersections
with the sets $\{trG_{s_{1},s_{2}}<R\}$ one can achieve that the
family $\Gamma$ will consist of the nuclear operators. Finally denote
by $\{e_{n};\: n\geq1\}$ the orthonormal basis in $H$ which is a
common eighenbasis for all operators from $\Gamma$ . Denote by $\check{\Gamma}$
the family of all projections onto subspaces generated by a finite
number of the vectors from $\{e_{n};\: n\geq1\}.$ Then every operator
$G_{s_{1},s_{2}}$takes values in $\check{\Gamma}.$ 

Now describe a random structure of $G_{s_{1},s_{2}}$. For the fixed
$n$ consider a random process in $H$ \[
\xi_{n}(t)=G_{0,t}e_{n},\: t\geq0.\]
By the definition of a random semigroup $\xi_{n}$ is a homogeneous
Markov process. From the other side there exists a random moment $\tau_{n}$
such that \[
\xi_{n}(t)=e_{n},\: t<\tau_{n},\:\xi_{n}(t)=0,\: t\geq\tau_{n}.\]
The random moment $\tau_{n}$ has an exponential distribution with
parameter $\lambda_{n}$. Define the random measure $\nu$ on the
product $[0;+\infty)\times\mathit{N}$ as follows

\[
\nu((s;t]\times A)=\sum_{n\in A}1_{\{G_{s,t}e_{n}=0\}}.\]
Note that the measure $\nu$ has independent values on the sets which
have disjoint projections on $[0;+\infty)$. For arbitrary $n\geq1$
the process $\{\nu((0;t]\times n);\: t\geq0\}$ is Poissonian with
the parameter $\lambda_{n}.$ Also note that for arbitrary $t>0$
\[
\mathit{P}\{\exists n_{0}\:\forall n\geq n_{0}:\:\xi_{n}(t)=0\}=1.\]
The theorem is proved. \textifsymbol[ifgeo]{48}

\section{{\normalsize Widths of compact sets defined with respect to the semigroups
of projections.}}

The last statement of the previous section gives us a description
of the semigroups of random finite-dimensional projections in terms
of integer-valued random measure. To understand the relationships
between this measure and geometrical properties of the semigroup consider
the asymptotic of widths of compact sets with respect to the images
of the semigroup projections. Let $\{G_{s,t},\:0\leq s\leq t<\infty\}$
be a random semigroup of finite-dimensional projections and $K$ be
a compact subset of $H$. We will investigate the behavior of the
value\[
\varsigma_{K}(t)=\max_{x\in K}||x-G_{0,t}x||\]
as $t\rightarrow0.$ The value $\varsigma_{K}(t)$ is exactly the
width of $K$ with respect to the linear subspace $G_{0,t}(H)$ {[}3].
Theorem 3.1 implies that with probability one $G_{0,t}$ strongly
converges to identity when $t\rightarrow0.$ Consequently, with probability
one $\varsigma_{K}(t)\rightarrow0,\: t\rightarrow0.$ We will investigate
the rate of the convergence. Let us consider the case, when the processes
$\{\xi_{n};\: n\geq1\}$ which arose in the description of the structure
of the semigroup are independent and the compact $K$ has a simple
description in the basis $\{e_{n};\: n\geq1\}.$ 

\textbf{\large Example 4.1.} Suppose that $\lambda_{n}=n,\: n\geq1$
and \[
K=\{x:\:(x,\: e_{n})^{2}\leq\frac{1}{n^{2}},\: n\geq1\}.\]
Now \[
\varsigma_{K}(t)^{2}=\sum_{n=1}^{\infty}\frac{\xi_{n}(t)}{n^{2}}.\]
 One can check that \[
\mathit{E}\varsigma_{K}(t)^{2}=\sum_{n=1}^{\infty}\frac{1}{n^{2}}(1-\exp(-nt))=\sum_{n=1}^{\infty}\int_{0}^{t}\int_{s}^{+\infty}\exp(-nr)drds=\]
 \[
\int_{0}^{t}\int_{s}^{+\infty}\sum_{n=1}^{\infty}\exp(-nr)drds=\int_{0}^{t}\int_{s}^{+\infty}\frac{\exp(-r)}{1-\exp(-r)}drds.\]
 Consequently, \[
\mathit{E}\varsigma_{K}(t)^{2}=\int_{0}^{t}ln(1-\exp(-s))ds\thicksim tlnt,\: t\rightarrow0.\]
In the similar way the forth moment can be estimated\[
\mathit{E}\varsigma_{K}(t)^{4}=\mathit{E}\sum_{i,j=1}^{\infty}\frac{\xi_{i}(t)\xi_{j}(t)}{i^{2}j^{2}}=\sum_{i=1}^{\infty}\frac{1}{i^{4}}(1-\exp(-it))+\]

\[
+\sum_{i\neq j}^{\infty}\frac{1}{i^{2}j^{2}}(1-\exp(-it))(1-\exp(-jt))=(\mathit{E}\varsigma_{K}(t)^{2})^{2}+\sum_{i=1}^{\infty}\frac{1}{i^{4}}(1-\exp(-it))\exp(-it).\]

Note that \[
\sum_{i=1}^{\infty}\frac{1}{i^{4}}(1-\exp(-it))\exp(-it)\leq\sum_{i=1}^{\infty}\frac{1}{i^{4}}(1-\exp(-it))\leq\sum_{i=1}^{\infty}\frac{t}{i^{3}}=ct.\]

Hence \[
\mathit{V}\varsigma_{K}(t)^{2}=o(tlnt),\: t\rightarrow0.\]

Finally, one can conclude that \[
\mathit{P}-\lim_{t\rightarrow0}\frac{\varsigma_{K}(t)}{\sqrt{tlnt}}=1.\]

In the next example we consider another type of compact and the same
semigroup.

\textbf{\large Example 4.2.} Consider the following compact set\[
K=\{x:\:\sum_{n=1}^{\infty}n^{2}(x,\: e_{n})^{2}\leq1\}.\]
For the same semigroup $\{G_{s,t},\:0\leq s\leq t<\infty\}$ as in
the previous example let us study the behavior of $\varsigma_{K}(t).$
Now\[
\varsigma_{K}(t)^{2}=\max_{n:\:\xi_{n}(t)=0}\frac{1}{n^{2}}.\]
Let us find the expectation \[
\mathit{E}\varsigma_{K}(t)^{2}=\sum_{n=1}^{\infty}\frac{1}{n^{2}}\prod_{i=1}^{n-1}\exp(-it)(1-\exp(-nt))=\sum_{n=2}^{\infty}\frac{1}{n^{2}}\exp(-\frac{n(n-1)}{2}t)\times\]
\[
\times(1-\exp(-nt))+1-\exp(-t).\]
The asymptotic behavior when $t\rightarrow0$ of the last summand
is trivial. To obtain the estimation from belaw let us rewrite the
sum as follows

\[
\sum_{n=2}^{\infty}\frac{1}{n^{2}}\exp(-\frac{n(n-1)}{2}t)(1-\exp(-nt))=\sum_{n=2}^{\infty}\frac{1}{n}\exp(-\frac{n(n-1)}{2}t)\int_{0}^{t}\exp(-ns)ds\geq\]
\[
\geq\sum_{n=2}^{\infty}\frac{1}{n}\exp(-\frac{n(n-1)}{2}t)t\exp(-nt)\geq t\int_{2}^{+\infty}\frac{1}{x}\exp(-\frac{x(x-1)}{2}t-xt)dx.\]
 For arbitrary sufficiently small positive $\alpha$ there exists
such positive $c$ that\[
\int_{c}^{+\infty}\frac{1}{x}\exp(-\frac{x(x-1)}{2}t-xt)dx\geq\int_{c}^{+\infty}\frac{1}{x}\exp(-(\frac{1}{2}+\alpha)x^{2}t)dx=\]
\[
=\int_{c\sqrt{(\frac{1}{2}+\alpha)t}}^{+\infty}\frac{1}{x}\exp(-x^{2})dx\thicksim\frac{1}{2}\ln t,\: t\rightarrow0.\]
For the upper estimate let us proceed in the same way\[
\sum_{n=2}^{\infty}\frac{1}{n^{2}}\exp(-\frac{n(n-1)}{2}t)(1-\exp(-nt))=\sum_{n=2}^{\infty}\frac{1}{n}\exp(-\frac{n(n-1)}{2}t)\int_{0}^{t}\exp(-ns)ds\leq\]
\[
\leq\sum_{n=2}^{\infty}\frac{1}{n}\exp(-\frac{n(n-1)}{2}t)t\leq t\int_{1}^{+\infty}\frac{1}{x}\exp(-\frac{x(x-1)}{2}t)dx\leq\]
\[
\leq t\int_{1}^{c}\frac{1}{x}\exp(-\frac{x(x-1)}{2}t)dx+t\int_{c}^{+\infty}\frac{1}{x}\exp(-(\frac{1}{2}-\alpha)x^{2}t)dx=\]

\[
=t\int_{1}^{c}\frac{1}{x}\exp(-\frac{x(x-1)}{2}t)dx+t\int_{c\sqrt{(\frac{1}{2}-\alpha)t}}^{+\infty}\frac{1}{x}\exp(-x^{2})dx\thicksim\frac{1}{2}t\ln t,\: t\rightarrow0.\]
To understand a piecewise behavior of $\varsigma_{K}(t)$ when $t\rightarrow0$
let us introduce a family $\{\tau_{n};\: n\geq1\}$ of independent
exponentially distributed random variables with intensities $n.$
Then the sequence of random processes $\{\xi_{n}(t);\: n\geq1\}$
is equidistributed with the sequence $\{1_{[0;\: t]}(\tau_{n});\: n\geq1\}.$
Consequently, for a continuous strictly decreasing positive function
$a$ and a constant $c>0$ \[
\mathit{P}\{\liminf_{t\rightarrow0}\frac{1}{\varsigma_{K}(t)a(t)}>c\}=\mathit{P}\{\exists N\:\forall n\geq N:\:\min_{j=1,\ldots,n}\tau_{j}>a^{-1}(\frac{n}{c})\}=\]
\[
\mathit{P}\{\exists N\:\forall n\geq N:\:\tau_{n}>a^{-1}(\frac{n}{c})\}.\]
 The last probability equals to one if and only if the infinite product
\[
\prod_{n=1}^{\infty}\mathit{P}\{\tau_{n}>a^{-1}(\frac{n}{c})\}\]
 converges. This condition is equivalent to the convergence of the
series\[
\sum_{n=1}^{\infty}\mathit{P}\{\tau_{n}<a^{-1}(\frac{n}{c})\}=\sum_{n=1}^{\infty}(1-\exp(-n)a^{-1}(\frac{n}{c})).\]
 This series converges if and only if \[
\sum_{n=1}^{\infty}na^{-1}(\frac{n}{c})<+\infty.\]
 This inequality holds simultaneously for all positive $c,$ which
means that\[
\lim_{t\rightarrow0}\varsigma_{K}(t)a(t)=0.\]
 For example, this condition is true for the function \[
a^{-1}(n)=\frac{1}{n^{2}\ln^{2}n}.\]
The upper bound can be obtained using the equality\[
\mathit{P}\{\frac{1}{\varsigma_{K}(t)}\geq n\}=\exp(-\frac{n(n-1)}{2}t),\: n\geq2.\]
 Define a function $\varphi$ by the formula\[
\varphi(t)=\sqrt{\frac{2}{t}llnt}\]
for sufficiently small positive $t$ with the usual agreement for
$lln.$ Then, taking $t_{n}=q^{n}$ for some $0<q<1$ one can get\[
\sum_{n=1}^{\infty}\mathit{P}\{\frac{1}{\varsigma_{K}(t_{n})}\geq(1+\varepsilon)\varphi(t_{n})\}<+\infty\]
 for any positive $\varepsilon.$ Consequently,\[
\limsup_{n\rightarrow\infty}\frac{1}{\varsigma_{K}(t_{n})\varphi(t_{n})}\leq1.\]
 Since the function $\varsigma_{K}$ is increasing, then for $t_{n+1}<t\leq t_{n}$
\[
\frac{1}{\varsigma_{K}(t)\varphi(t)}\leq\frac{1}{\varsigma_{K}(t_{n+1})\varphi(t_{n+1})}\frac{\varphi(t_{n+1})}{\varphi(t_{n})}.\]
Hence\[
\limsup_{t\rightarrow0}\frac{1}{\varsigma_{K}(t)\varphi(t)}\leq\frac{1}{q}.\]
Finally\[
\liminf_{t\rightarrow0}\varsigma_{K}(t)\varphi(t)\geq1.\]

In general the structure of $K$ can be more complicated and does
not allow an explicit form for $\varsigma_{K}(t).$ In some cases
one can have only the estimation for the Kolmogorov width for $K$
\[
d_{n}(K)=\inf_{dimL=n}\max_{x\in K}\rho(x,\: K),\]
where $\inf$ is taken over all subspaces $L$ of $H,$ which have
the dimension $n.$ From this reason it is useful to estimate the
growth of $\dim G_{0,t}(H)$ when $t\rightarrow0.$ For the semigroup
from the previous examples such an estimation can be obtained as follows.

\textbf{\large Example 4.3.} Define $\alpha(t)=\dim G_{0,t}(H).$
Using the random variables $\{\tau_{n};\: n\geq1\}$ which were defined
in Example 4.2 one can check, that\[
\mathit{P}\{\alpha(t)\geq n\}=\mathit{P}\{\exists\: k_{1}<k_{2}<\ldots<k_{n}:\:\tau_{k_{1}}\geq t,\:\tau_{k_{2}}\geq t,\ldots,\:\tau_{k_{n}}\geq t\}.\]
 Consequently, \[
\mathit{P}\{\alpha(t)\geq n\}\leq\mathit{P}\{\tau_{1}\geq t,\:\tau_{2}\geq t,\ldots,\:\tau_{n}\geq t\}+1-\mathit{P}\{\forall\: k>n:\:\tau_{k}<t\}\leq\]
\[
\leq\exp(-nt)+1-\prod_{k=n+1}^{\infty}(1-\exp(-kt))\leq\exp(-nt)+\sum_{k=n+1}^{\infty}\exp(-kt)\leq\]
\[
\leq\exp(-nt)+\frac{\exp(-nt)}{1-\exp(-t)}\leq\exp(-nt)(1+(1-e^{-1})\frac{1}{t})\]
 for $t\in(0;\: t).$ Taking $t_{k}=\frac{1}{k},\: n_{k}=[(2+\delta)k\ln k]$
for positive $\delta$ one can get that\[
\limsup_{k\rightarrow\infty}\frac{\alpha(\frac{1}{k})}{k\ln k}\leq2.\]
Here $[x]$ means an integer part of $x.$ Using the monotonicity
of $\alpha$ one can conclude that with probability one\[
\limsup_{t\rightarrow0}\frac{t\alpha(t)}{2|\ln t|}\leq1.\]
 To obtain an estimation from below let us denote \[
c(t)=\prod_{j=1}^{\infty}(1-\exp(-jt))=\mathit{P}\{\alpha(t)=0\}.\]
 Note, that\[
\ln c(t)=\sum_{j=1}^{\infty}\ln(1-\exp(-jt))\leq-\sum_{j=1}^{\infty}\exp(-jt)=-\frac{\exp(-t)}{1-\exp(-t)}\sim-\frac{1}{t},\: t\rightarrow0.\]
 Consequently, \[
\lim_{t\rightarrow0}t\ln c(t)=-1.\]
Now

\[
\mathit{P}\{\alpha(t)<n\}=\sum_{k=0}^{n-1}\mathit{P}\{\alpha(t)=k\}=\prod_{j=1}^{\infty}(1-\exp(-jt))(1+\]
 \[
+\sum_{k=1}^{n-1}\sum_{1\leq j_{1}<j_{2}<\ldots<j_{k}}\prod_{s=1}^{k}\exp(-j_{s}t)(1-\exp(-j_{s}t))^{-1}).\]
 Consider the series\[
\sum_{j=2}^{\infty}\exp(-jt)(1-\exp(-jt))^{-1}\leq\int_{1}^{\infty}\exp(-xt)(1-\exp(-xt))^{-1}dx=\]
\[
=-\frac{1}{t}\ln(1-\exp(-t)).\]
 Hence\[
\mathit{P}\{\alpha(t)<n\}\leq\prod_{j=1}^{\infty}(1-\exp(-jt))\sum_{k=0}^{n-1}\frac{1}{k!}(\exp(-t)(1-\exp(-t))^{-1}-\frac{1}{t}\ln(1-\exp(-t)))^{k}.\]
 Consequently, for arbitrary $\delta>1,\: c>1$ for sufficiently small
$t$ \[
\mathit{P}\{\alpha(t)<n\}\leq c(t)c^{n}\frac{1}{t^{n}}|\ln t|^{n}\leq\]
\[
\leq\exp(-\frac{\delta}{t}+n\ln c+n|\ln t|+n\ln|\ln t|).\]
For $\beta>2$ consider the sequences $\{t_{k}=\frac{1}{k};\: k\geq1\}$
and $\{n_{k}=\frac{k}{\beta\ln k};\: k\geq2\}.$ Then \[
\sum_{k=2}^{\infty}\mathit{P}\{\alpha(t_{k})<n_{k}\}<+\infty.\]
 Using the monotonicity of $\alpha$ as above one can get that with
probability one\[
\liminf_{t\rightarrow0}\alpha(t)t|\ln t|\geq\frac{1}{2}.\]

\end{document}